\input amstex
\input amsppt.sty
\magnification=\magstep1
\hsize=30truecc
\vsize=22.2truecm
\baselineskip=16truept
\TagsOnRight
\pageno=1
\nologo
\def\Z{\Bbb Z}
\def\N{\Bbb N}

\def\Q{\Bbb Q}

\def\l{\left}
\def\r{\right}
\def\bg{\bigg}
\def\({\bg(}
\def\[{\bg\lfloor}
\def\){\bg)}
\def\]{\bg\rfloor}
\def\t{\text}
\def\f{\frac}

\def\se {\subseteq}

\def\bi{\binom}
\def\eq{\equiv}

\def\ls{\leqslant}

\def\mo{\roman{mod}}
\def\ord{\roman{ord}}

\def\ve{\varepsilon}
\def\al{\alpha}

\def\Proof{\noindent{\it Proof}}

\def\Remark{\medskip\noindent{\it  Remark}}

\def\Ack{\medskip\noindent {\bf Acknowledgment}}
\hbox {Proc. Amer. Math. Soc. 139(2011), no.\,5, 1569--1577.}
\bigskip
\topmatter
\title Binomial coefficients and the ring of $p$-adic integers\endtitle
\author Zhi-Wei Sun* and Wei Zhang\endauthor
\leftheadtext{Zhi-Wei Sun and Wei Zhang}
\affil Department of Mathematics, Nanjing University\\
 Nanjing 210093, People's Republic of China
 \\ {\tt zwsun$\@$nju.edu.cn, \ \ zhangwei$\_$07$\@$yahoo.com.cn}
\endaffil
\abstract Let $k>1$ be an integer and let $p$ be a prime. We show
that if  $p^a\ls k<2p^a$ or $k=p^aq+1$ (with $q<p/2$) for some
$a=1,2,3,\ldots$, then the set $\{\bi nk:\, n=0,1,2,\ldots\}$ is
dense in the ring $\Z_p$ of $p$-adic integers, i.e., it contains a
complete system of residues modulo any power of $p$.
\endabstract
\thanks 2010 {\it Mathematics Subject Classification}.\,Primary 11B65;
Secondary 05A10, 11A07, 11S99.
\newline\indent *This author is the corresponding author. He is supported
by the National Natural Science Foundation (grant 10871087) and the
Overseas Cooperation Fund (grant 10928101) of China.
\endthanks
\endtopmatter
\document

\heading{1. Introduction}\endheading

Let $p$ be an odd prime. In section F11 of Guy [Gu, p.\,381] it is conjectured
that $|\{n!\ \mo\ p:\ n=1,2,3,\ldots\}|$ is about $p(1-1/e)$
asympototically. [CVZ] provided certain evidence for the conjecture.

In [BLSS] the authors proved that for infinitely many primes $p$
there are at least $\log\log p/\log\log\log p$ distinct integers
among $0,1,\ldots,p-1$ which are not congruent to $n!$ for any
$n\in\Z^+=\{1,2,3,\ldots\}$.

Garaev and Luca [GL] showed that for any $\ve>0$ there is a
computable positive constant $p_0(\ve)$ such that  for any prime
$p>p_0(\ve)$ and integers $t> p^{\ve}$ and
$s>t+p^{1/4+\ve}$ we have
$$\{m_1!\cdots m_t!\ (\mo\ p):\ m_1+\cdots+m_t=s\}\supseteq\{r\ (\mo\
p):\ r=1,\ldots,p-1\}.$$

Let $p$ be any prime. As usual, we denote by $\Z_p$ the ring of
$p$-adic integers in the $p$-adic field $\Q_p$. The reader may
consult an excellent book [M] by Murty for the basic knowledge
of $p$-adic analysis. Any given $p$-adic integer $\al$ has a unique
representation in the form
$$\al=\sum_{j=0}^{\infty}a_jp^j\ \ \t{with}\ a_j\in[0,p-1]=\{0,1,\ldots,p-1\}.$$
For each $b\in\N=\{0,1,2,\ldots\}$, we have
$$\al\eq r(b)\ (\mo\ p^b),\ \ \t{i.e.},\ \ |\al-r(b)|_p\ls\f1{p^b},$$
where $r(b):=\sum_{0\ls j<b}a_jp^j$ and $|\cdot|_p$ is the $p$-adic
norm.

In this paper we study the following new problem (which was actually
motivated by the first author's paper [S]).

\proclaim{Problem 1.1} Given a prime $p$ and a positive integer $k$,
is the set $\{\bi nk:\ n\in\N\}$ dense in $\Z_p$? In other
words, does the set contain a complete system of residues modulo any
power of $p$?
\endproclaim

\medskip
\noindent{\bf Definition 1.1.} Let $k\in\N$ and
$m\in\Z^+$, and define
$$R_m(k):=\bg\{\bi nk\ (\mo\ m):\, n\in\N\bg\}.\tag1.1$$
If $R_m(k)=\Z/m\Z$, then we call $m$ a {\it $k$-universal} number.
\medskip

Clearly all positive integers are $1$-universal and $1$ is $k$-universal for all $k\in\Z^+$.

If $p$ is a prime, $a,k,n,n'\in\N$ and $n'\eq n\ (\mo\
p^{a+\ord_p(k!)})$ then
$$\bi {n'}k=\f{\prod_{0\ls j<k}(n'-j)}{k!}\eq\f{\prod_{0\ls j<k}(n-j)}{k!}
=\bi nk\ (\mo\ p^a).$$ Combining this observation with the Chinese
Remainder Theorem we immediately get the following basic
proposition.

\proclaim{Proposition 1.1} Let $k\in\Z^+$ and $m=p_1^{a_1}\cdots
p_r^{a_r}$, where $p_1,\ldots,p_r$ are distinct primes and
$a_1,\ldots,a_r\in\Z^+$. Then $m$ is $k$-universal if and only if
$p_1^{a_1},\ldots,p_r^{a_r}$ are all $k$-universal.
\endproclaim

In view of Proposition 1.1, we may focus on those $k$-universal prime powers.

Let $k>1$ be an integer. If $p>k$ is a prime, then $p$ is not $k$-universal since
$$\bg\{\bi 0k,\bi 1k,\ldots,\bi{p-1}k\bg\}$$
is not a complete system of residues modulo $p$. (Note that $\bi 0k=\bi1k=0$.)
Thus, if $m\in\Z^+$ is $k$-universal, then
$m$ has no prime divisor greater than $k$.

For an integer $k>1$, a prime $p>k$ and an integer $r\in[1,p-1]=\{1,\ldots,p-1\}$,
the congruence $\bi xk\eq r\ (\mo\ p)$ might have
more than two solutions. For example,
$$\bi{12}5\eq\bi{19}5\eq\bi{22}5\eq\bi{31}5\eq18\ (\mo\ 43)$$
and $$\bi{15}{10}\eq\bi{21}{10}\eq\bi{25}{10}\eq\bi{30}{10}\eq14\ (\mo\ 61).$$

Recall the following useful result of Lucas.

\proclaim{Lucas' Theorem {\rm (cf. [Gr] and [HS])}} Let $p$ be any
prime, and let $n_0,k_0,\ldots,n_r,k_r\in[0,p-1]$. Then we have
$$\bi{\sum_{i=0}^rn_ip^i}{\sum_{i=0}^rk_ip^i}\eq\prod_{i=0}^r\bi{n_i}{k_i}\
(\mo\ p).$$
\endproclaim

Clearly Lucas' theorem implies the following proposition.
\proclaim{Proposition 1.2} Let $p$ be a prime and let $k=\sum_{i=0}^r k_ip^i$ with $k_i\in[0,p-1]$.
Then $R_p(k)=\prod_{i=1}^rR_p(k_i)$. In particular, when $k_0,\ldots,k_r\in\{0,p-2,p-1\}$
we have
$$R_p(k)\se\{r(\mo\ p):\ r=0,\pm1\}.$$
\endproclaim

Let $p$ be a prime. As $R_p(1)=\Z/p\Z$, if the $p$-adic expansion of $k\in\Z^+$ has a digit 1 then
$R_p(k)=\Z/p\Z$ by Proposition 1.2. In this spirit, Proposition 1.2 is helpful to study when $R_p(k)=\Z/p\Z$.
For an integer $b>1$, to investigate whether $p^b$ is $k$-universal (i.e., $R_{p^b}(k)=\Z/p^b\Z$) one might think that
we should use extended Lucas theorem for prime powers.
However, all known generalizations of Lucas' theorem to prime powers
are somewhat unnatural and complicated, e.g., K. Davis and W. Webb [DW] proved that if $p>3$ is a prime,
$a,b,n,k\in\Z^+$ and $n_0,k_0\in[0,p^a-1]$ then
$$\bi{np^{a+b}+n_0}{kp^{a+b}+k_0}\eq\bi{np^{\lfloor b/3\rfloor}}{kp^{\lfloor b/3\rfloor}}\bi{n_0}{k_0}\ (\mo\ p^{b+1}).$$
Therefore, we prefer to approach Problem 1.1 by induction argument which can be easily understood.

Our first result is as follows.

\proclaim{Theorem 1.1} Let $p$ be a prime and let
$a\in\N=\{0,1,2,\ldots\}$. Let $k$ be an integer with $p^a\ls
k<2p^a$. Then, for any $b\in\N$ and $r\in\Z$ there is an integer
$n\in[0,p^{a+b}-1]$ with $n\eq k\ (\mo\ p^{a})$ such that $\bi nk\eq
r\ (\mo\ p^b)$.
\endproclaim
\Remark\ 1.1. For a prime $p$ and a positive integer $k$ having a digit 1 in its $p$-adic expansion,
$p$ is definitely $k$-universal (i.e., $R_p(k)=\Z/p\Z$) but
$p^b$ might be not  for some integer $b>1$. For example,
$21=4\times5+1$ and $\ord_5(21!)=4$, thus
$$\align\l\{\bi n{21}\ (\mo\ 5^2):\ n\in\N\r\}=&\l\{\bi n{21}\ (\mo\ 5^2):\ n\in[0,5^6-1]\r\}
\\=&\{r\, (\mo\ 5^2):\ r=0,\pm1,\pm3,\pm5,\pm10\}
\endalign$$
and hence $R_{5^2}(21)\not=\Z/5^2\Z$.

\medskip

 Here is a consequence of
Theorem 1.1.

\proclaim{Corollary 1.1}  Let $p$ be a prime and let $k\in\Z^+$ with
$\log_p(k/2)<\lfloor\log_pk\rfloor$. Then $1,p,p^2,\ldots$ are
$k$-universal numbers and the set $\{\bi nk:\,n\in\N\}$ is a dense subset
of the ring $\Z_p$ of $p$-adic integers.
\endproclaim
\Proof. Set $a=\lfloor\log_pk\rfloor$. Then $p^a\ls k<2p^a$. By Theorem 1.1, $p^b$ is $k$-universal
for every $b=0,1,2,\ldots$. Therefore $\{\bi nk:\, n\in\N\}$ is dense in $\Z_p$.
This concludes the proof. \qed

For any $k\in\Z^+$ there is a unique $a\in\N$ such that $2^a\ls k<2^{a+1}$.
Thus Theorem 1.1 or Corollary 1.1 implies the following result.

\proclaim{Corollary 1.2} Let $k\in\Z^+$. Then any power of two is $k$-universal
and hence the set $\{\bi nk:\,n\in\N\}$ is a dense subset of the $2$-adic integral ring $\Z_2$.
\endproclaim

\medskip\noindent{\bf Definition 1.2.} A positive integer $k$ is said to be {\it universal}
if any power of a prime $p\ls k$ is $k$-universal, i.e., $\{\bi
nk:\,n\in\N\}$ is a dense subset of $\Z_p$ for any prime $p\ls k$.
\medskip

Theorem 1.1 implies that 1, 2, 3, 4, 5, 9 are universal numbers. To obtain other universal numbers, we need to extend Theorem 1.1.

\proclaim{Theorem 1.2} Let $p$ be a prime and let $a\in\N$. Let
$k=k_0+p^ak_1$ with $k_0\in[0,p^a-1]$ and $k_1\in[1,p-1]$. Suppose
that for each $r=1,\ldots,p-1$ there are $n_0\in[k_0,p^a-1]$ and
$n_1\in[k_1,p-1]$ such that
 $$\bi {n_1}{k_1}\bi{n_0}{k_0}\eq r\ (\mo\ p)\ \ \t{and}\ \ P_{k_1}(n_1)\not\eq0\ (\mo\ p),\tag 1.2$$
where
 $$P_{k_1}(x)=\sum_{j=1}^{k_1}\f{(-1)^{j-1}}j\bi{x}{k_1-j}.\tag1.3$$
Then, for any $b\in\N$, the set $\{\bi nk:\, n\in[0,p^{a+b}-1]\}$ contains a complete system of residues modulo $p^b$.
\endproclaim

\Remark\ 1.2. Let $p$ be an odd prime. Then
$$\align P_{p-1}(p-1)=&\sum_{j=1}^{p-1}\f{(-1)^{j-1}}j\bi{p-1}j
\\\eq&-\sum_{j=1}^{p-1}\f1j=-\sum_{j=1}^{(p-1)/2}\l(\f1j+\f1{p-j}\r)\eq0\ (\mo\ p).
\endalign$$
Thus, for $k_1=p-1$ there is no $n_1\in[k_1,p-1]$ with $P_{k_1}(n_1)\not\eq0\ (\mo\ p)$.

\proclaim{Corollary 1.3} Let $p$ be an odd prime and $q\in\{1,\ldots,(p-1)/2\}$.
Then, for any $a\in\Z^+$ and $b\in\N$, the number $p^b$ is $(p^aq+1)$-universal.
\endproclaim
\Proof. Let $k_1=q$, $k_0=1$, and $k=p^ak_1+k_0=p^aq+1$. As $P_{k_1}(x)\eq0\ (\mo\ p)$
cannot have more than $\deg P_{k_1}(x)=k_1-1$ solutions (see, e.g., [IR, p.\,39])
there exists $n_1\in[k_1,2k_1-1]\se[k_1,p-1]$
such that $P_{k_1}(n_1)\not\eq0\ (\mo\ p)$. Note that $\bi{n_1}{k_1}\not\eq0\ (\mo\ p)$.
For any $r\in[1,p-1]$ there is a unique $n_0\in[1,p-1]$ such that
$$\bi{n_1}{k_1}\bi{n_0}{k_0}=n_0\bi{n_1}{k_1}\eq r\ (\mo\ p).$$
Applying Theorem 1.2, we immediately obtain the desired result. \qed

From Theorem 1.2 we can deduce the following result.

\proclaim{Theorem 1.3} The integers $11$, $17$ and $29$ are universal
numbers.
\endproclaim

We have the following conjecture based on our computation via the software {\tt Mathematica}.

\proclaim{Conjecture 1.1} There are no universal numbers other than {\rm 1, 2, 3, 4, 5, 9, 11, 17, 29}.
\endproclaim

In Sections 2, 3 and 4 we will prove Theorems 1.1, 1.2 and 1.3 respectively.

\heading{2. Proof of Theorem 1.1}\endheading

 We use induction on $b$. The case $b=0$ is trivial, so we proceed to the induction step.

 Now fix $b\in\N$ and $r\in\Z$. Suppose that $m\in\Z$, $n=k+p^am\in[0,p^{a+b}-1]$ and
 $\bi nk\eq r\ (\mo\ p^b)$. Let $q$ be the smallest nonnegative residue of
 $(r-\bi nk)/p^b$ modulo $p$.

 Set
 $n'=n+p^{a+b}q$. Then
 $$n'<p^{a+b}(q+1)\ls p^{a+b+1}\ \ \t{and}\ \ n'\eq n\eq k\ (\mo\ p^a).$$
 By the Chu-Vandermonde identity (cf. (5.22) of [GKP, p.\,169]),
 $$\bi {n'}k=\bi{n+p^{a+b}q}k=\sum_{j=0}^k\bi{p^{a+b}q}j\bi
 n{k-j}.$$
 If $1\ls j\ls k$ and $j\not=p^a$ then $p^a\nmid j$ and hence
 $$\bi{p^{a+b}q}j=\f{p^{a+b}q}j\bi{p^{a+b}q-1}{j-1}\eq0\ (\mo\
 p^{b+1}).$$ Note also that
 $$\bi{p^{a+b}q}{p^a}=p^bq\prod_{t=1}^{p^a-1}\f{p^{a+b}q-t}t\eq
 p^bq(-1)^{p^a-1}\eq p^bq\ (\mo\ p^{b+1}).$$
 Therefore
 $$\bi{n'}k\eq\bi nk+ p^bq\bi  n{k-p^a}\eq r-p^bq+p^bq\bi  n{k-p^a}\ (\mo\ p^{b+1}).$$
So it suffices to show that
$$\bi n{k-p^a}\eq1\ (\mo\ p).$$

 As $0\ls k-p^a<p^a$, Lucas' theorem implies that
$$\bi n{k-p^a}=\bi{(m+1)p^a+(k-p^a)}{0p^a+(k-p^a)}\eq\bi{m+1}0\bi{k-p^a}{k-p^a}=1\ (\mo\ p).$$

 Combining the above we have completed the proof by induction.

\heading{3. Proof of Theorem 1.2}\endheading

We claim that for each $b=0,1,2,\ldots$ the set $\{\bi nk:\, n\in S(b)\}$
contains a complete system of residues modulo $p^b$, where
$$S(b)=\bg\{n\in[0,p^{a+b}-1]:\ \bi{\{n\}_{p^a}}{k_0}
\sum_{j=1}^{k_1}\f{(-1)^{j-1}}j\bi{\lfloor n/p^a\rfloor}{k_1-j}\not\eq0\ (\mo\ p)\bg\}$$
and $\{n\}_{p^a}$ denotes the least nonnegative residue of $n$ mod $p^a$.

The claim is trivial for $b=0$ since
$$\bi{\{k_0\}_{p^a}}{k_0}\sum_{j=1}^{k_1}\f{(-1)^{j-1}}j
\bi{\lfloor k_0/p^a\rfloor}{k_1-j}=\f{(-1)^{k_1-1}}{k_1}\not\eq0\ (\mo\ p).$$

As $\deg P_{k_1}(x)<k_1$, there exists  $n_1\in[0,k_1-1]$ such that
$P_{k_1}(n_1)\not\eq0\ (\mo\ p)$. Combining this with the
supposition in Theorem 1.2, we see that for any $r\in[0,p-1]$ there
are $n_0\in[0,p^a-1]$ and $n_1\in[0,p-1]$ satisfying (1.2) and the
congruence $\bi{n_0}{k_0}\not\eq0\ (\mo\ p)$. Taking
$n=p^an_1+n_0\in[0,p^{a+1}-1]$ we find that
$$\bi nk\eq\bi {n_1}{k_1}\bi{n_0}{k_0}\eq r\ (\mo\ p)$$
by Lucas' theorem. This proves the claim for $b=1$.

Now let $b\in\Z^+$ and assume that $\{\bi nk:\, n\in S(b)\}$ contains a complete system of residues modulo $p^b$.
We proceed to prove the claim for $b+1$.

Let $r$ be any integer. By the induction hypothesis, there is an
integer $n\in[0,p^{a+b}-1]$ such that
$$\bi nk\eq r\ (\mo\ p^b)\ \ \t{and}\ \
\bi{n_0}{k_0}\sum_{j=1}^{k_1}\f{(-1)^{j-1}}j\bi{\lfloor
n/p^a\rfloor}{k_1-j}\not\eq0\ (\mo\ p),$$ where $n_0=\{n\}_{p^a}$.
Hence, for some  $q\in[0,p-1]$ we have
$$q\bi{n_0}{k_0}\sum_{j=1}^{k_1}\f{(-1)^{j-1}}j\bi{\lfloor n/p^a\rfloor}{k_1-j}\eq\f{r-\bi nk}{p^b}\ (\mo\ p).$$
Clearly, $n'=n+p^{a+b}q\in[0,p^{a+b+1}-1]$ and
$$\align&\bi{\{n'\}_{p^a}}{k_0}\sum_{j=1}^{k_1}\f{(-1)^{j-1}}j\bi{\lfloor n'/p^a\rfloor}{k_1-j}
\\=&\bi{n_0}{k_0}\sum_{j=1}^{k_1}\f{(-1)^{j-1}}j\bi{\lfloor n/p^a\rfloor+p^bq}{k_1-j}
\\\eq&\bi{n_0}{k_0}\sum_{j=1}^{k_1}\f{(-1)^{j-1}}j\bi{\lfloor n/p^a\rfloor}{k_1-j}
\not\eq0\ (\mo\ p).
\endalign$$
As in the proof of Theorem 1.1, we have
$$\align\bi {n'}k-\bi nk=&\sum_{j=1}^k\bi{p^{a+b}q}j\bi n{k-j}
\\\eq&\sum_{j=1}^{\lfloor k/p^a\rfloor}\bi{p^{a+b}q}{p^aj}\bi n{k-p^aj}\ (\mo\ p^{b+1}).
\endalign$$
By Lucas' theorem, for $1\ls j\ls \lfloor k/p^a\rfloor=k_1$ we have
$$\bi{p^{a+b}q}{p^aj}\eq\bi{p^bq}j=\f{p^bq}j\prod_{0<i<j}\f{p^bq-i}i\eq p^bq\f{(-1)^{j-1}}j\ (\mo\ p^{b+1})$$
and
$$\bi{n}{k-p^aj}=\bi{p^a\lfloor n/p^a\rfloor+n_0}{p^a(k_1-j)+k_0}
\eq\bi{\lfloor n/p^a\rfloor}{k_1-j}\bi{n_0}{k_0}\ (\mo\ p).$$
Therefore
$$\bi{n'}k-\bi nk\eq p^bq\bi{n_0}{k_0}\sum_{j=1}^{k_1}\f{(-1)^{j-1}}j\bi{\lfloor n/p^a\rfloor}{k_1-j}
\eq r-\bi nk\ (\mo\ p^{b+1})$$
and hence $\bi{n'}k\eq r\ (\mo\ p^{b+1})$. This concludes the induction step.

 In view of the above we have proved the claim and hence the desired result follows.

\heading{4. Proof of Theorem 1.3}\endheading

(I) We first prove that $11$ is universal.

Since
$$2^3<11<2^4,\ \ 3^2<11<2\times 3^2,\ \ 7<11<2\times7,$$
and $11=2\times 5+1$ with $2\ls(5-1)/2$, by Theorem 1.1 and
Corollary 1.3, 11 is universal.

\medskip

(II) Now we want to show that $17$ is universal.

Observe that
$$2^4<17<2^5,\ \ 3^2<17<2\times 3^2,\ \ 11<17<2\times11,$$
and $13<17<3\times 13$.
By Theorem 1.1, $p^b$ is $17$-universal for any $p=2,3,11,13$ and $b\in\N$.

Note that
$$\sum_{j=1}^{\lfloor 17/5\rfloor}\f{(-1)^{j-1}}j\bi{x}{\lfloor 17/5\rfloor-j}=\f{x^2-2x}2+\f13
\eq\f{(x-1)^2-2}2\not\eq0\ (\mo\ 5).$$
Also, $17=3\times5+2$, and
$$\align&\bi 33\bi22\eq1\ (\mo\ 5),\ \ \bi 33\bi32\eq3\ (\mo\ 5),
\\&\bi 43\bi22\eq4\ (\mo\ 5),\ \ \bi 43\bi32\eq2\ (\mo\ 5).
\endalign$$
So, $5^b$ is also 17-universal for any $b\in\N$.

Clearly
$$\sum_{j=1}^{\lfloor 17/7\rfloor}\f{(-1)^{j-1}}j\bi{x}{\lfloor 17/7\rfloor-j}=x-\f12\eq x-4\ (\mo\ 7).$$
Also, $17=2\times7+3$, and
$$\align&\bi 22\bi33\eq1\ (\mo\ 7),\ \ \bi 22\bi43\eq4\ (\mo\ 7),\ \ \bi 22\bi53\eq3\ (\mo\ 7),
\\&\bi 22\bi63\eq6\ (\mo\ 7),\ \ \bi 32\bi43\eq5\ (\mo\ 7),\ \ \bi 32\bi53\eq2\ (\mo\ 7).
\endalign$$
Thus, $7^b$ is also 17-universal for any $b\in\N$.

\medskip

(III) Finally we prove that 29 is universal.

By Theorem 1.1, it remains to prove that $p^b$ is 29-universal for any $p=7,11,13$ and $b\in\N$.

Note that $29=4\times7+1$. It is easy to check that
$$\sum_{j=1}^4\f{(-1)^{j-1}}j\bi{4}{4-j}\not\eq0\ (\mo\ 7).$$
For any $r\in[1,6]$, we have $\bi{4}{4}\bi{r}1\eq r\ (\mo\ 7)$.
So, by Theorem 1.2, $7^b$ is 29-universal for any $b\in\N$.

Clearly $29=2\times11+7$, and
$$\sum_{j=1}^2\f{(-1)^{j-1}}j\bi x{2-j}= x-\f12\eq x-6\ (\mo\ 11).$$
Observe that
$$\align&\bi 22\bi77\eq1\ (\mo\ 11),\ \bi 22\bi87\eq-3\ (\mo\ 11),
\\& \bi 22\bi97\eq3\ (\mo\ 11),\ \bi 22\bi{10}7\eq-1\ (\mo\ 11),
\\& \bi 32\bi87\eq2\ (\mo\ 11),\  \bi 32\bi97\eq-2\ (\mo\ 11),
\\&\bi 42\bi77\eq-5\ (\mo\ 11),\  \bi 42\bi{10}7\eq5\ (\mo\ 11),
\\&\bi 42\bi87\eq4\ (\mo\ 11),\ \bi 42\bi97\eq-4\ (\mo\ 11).
\endalign$$
Applying Theorem 1.2 we see that $11^b$ is 29-universal for any $b\in\N$.

Observe that $29=2\times 13+3$ and
$$\sum_{j=1}^2\f{(-1)^{j-1}}j\bi{x}{2-j}=x-\f12\eq x-7\ (\mo\ 13).$$
Also,
$$\align&\bi 22\bi33\eq1\ (\mo\ 13),\ \bi 22\bi43\eq4\ (\mo\ 13),
\\& \bi 22\bi53\eq-3\ (\mo\ 13),\ \bi 22\bi{6}3\eq-6\ (\mo\ 13),
\\&\bi22\bi 73\eq-4\ (\mo\ 13),\ \bi 22\bi93\eq6\ (\mo\ 13),
\\&\bi 22\bi{10}3\eq3\ (\mo\ 13),\  \bi 22\bi{12}3\eq-1\ (\mo\ 13),
\\&\bi 32\bi63\eq-5\ (\mo\ 13),\ \bi 32\bi93\eq5\ (\mo\ 13),
\\&\bi42\bi 43\eq-2\ (\mo\ 13),\ \bi{4}2\bi 73\eq2\ (\mo\ 13).
\endalign$$
Thus, with the help of Theorem 1.2, $13^b$ is 29-universal for any $b\in\N$.

\medskip

By the above, we have completed the proof of Theorem 1.3.

\medskip

\Ack. The authors are grateful to the referee for many helpful comments.

\enddocument